\documentclass[11pt]{article}
\usepackage{amsfonts}
\usepackage{latexsym}
\usepackage{amscd,mathrsfs,amsfonts}
\usepackage{amssymb,amsmath}
\usepackage{stmaryrd}
\usepackage{graphicx}

\newtheorem{theorem}{\bf Theorem}
\newtheorem{lemma}{\bf Lemma}

\newtheorem{remark}{\bf Remark}
\newtheorem{corollary}{\bf Corollary}
\newtheorem{proposition}{\bf Proposition}
\newtheorem{definition}{{\bf Definition}}

\begin{document}

\title{Invariant Measures with Bounded Variation Densities for Piecewise Area Preserving Maps}
\author{Yiwei Zhang and Congping Lin\\ \small {Mathematics Research Institute, University of Exeter, Exeter, EX4 4QF, UK}\\
(e-mail: yz297@exeter.ac.uk, cl336@exeter.ac.uk)}
\date{}
\maketitle
%

\begin{abstract}
We investigate the properties of absolutely continuous invariant
probability measures (ACIPs), especially those measures with bounded
variation densities, for piecewise area preserving maps (PAPs) on
$\mathbb{R}^d$. This class of maps unifies piecewise isometries
(PWIs) and piecewise hyperbolic maps where Lebesgue measure is
locally preserved. Using a functional analytic approach, we first
explore the relationship between topological transitivity and
uniqueness of ACIPs, and then give an approach to construct
invariant measures with bounded variation densities for PWIs. Our
results ``partially'' answer one of the fundamental questions posed
in \cite{Goetz03} - to determine all invariant non-atomic
probability Borel measures in piecewise rotations. When restricting
PAPs to interval exchange transformations (IETs), our results imply
that for non-uniquely ergodic IETs with two or more ACIPs, these
ACIPs have very irregular densities, i.e., they have unbounded
variation.

\end{abstract}


\section{Introduction}
Conservative systems are often used as models of the physical world,
where conservative is usually understood as energy preserving (i.e.,
where energy is invariant under the time evolution). In this article
we consider conservative systems that are governed by discrete time
dynamical systems. In particular, we focus on multidimensional
{\em piecewise area preserving maps} (PAPs), which is a general extension
of interval exchange transformations (IETs) into $\mathbb{R}^{d}$.
Regarding IETs, Keane conjectured that minimality implies unique
ergodicity in \cite{Keane75} and this conjecture holds for IETs with
two or three intervals. However, counterexamples have been
constructed; see \cite{Keane77,Keynes76}. Thereafter, Masur
\cite{Masur82} and Veech \cite{Veech82} have independently demonstrated
that almost every topologically transitive IET (with respect to Lebesgue
measure) is uniquely ergodic; simultaneously, Keane $\&$ Rauzy
\cite{Keane82} revealed that unique ergodicity holds for a Baire
residual subset of the space of IETs. To fully understand the
densities of absolutely continuous invariant measures (ACIPs) for these non-uniquely ergodic counterexamples, it is natural to explore equivalent conditions to the topological transitivity in IETs in terms of ACIPs.

When extending to the multidimensional PAPs, we are facing at least
two technical obstacles: complicated topology in high dimensions and
non-local preservation of distance. For the class of PAPs which
preserve distance locally, a special case of interest is the class
of piecewise isometries (PWIs). Establishing the properties of their
ACIPs will partially contribute to answering a fundamental question
posed in \cite{Goetz03}, i.e., to determine all invariant non-atomic
probability Borel measures for piecewise rotations. This question is
still open so far.

For the class of PAPs that do not preserve distance locally, a
particular case is piecewise hyperbolic maps. For these maps,
properties that have been studied include transitivity and
possession of a unique physical measure (e.g., see works of Boyarsky
$\&$ G\'{o}ra \cite{Boy97} and Viana \cite{Viana97}). These studies
use a functional analytic approach by choosing a ``reasonable''
function space and applying a {\em transfer operator} on this space.
They study statistical properties of the system by looking at the
operator fixed point and determining if there is a spectral gap. In
one-dimensional piecewise expanding maps, the space of bounded
variation functions has been demonstrated to be such a
``reasonable'' space \cite{Boy97,Viana97}. In higher dimensions, the
space of multidimensional bounded variation functions can still be
chosen under certain assumptions \cite{Buzzi01,Keller96,Tsujii00}
and contains a classical anisotropic Sobolev space of
Triebel-Lizorkin type \cite{Baladi09}.

In this article, our interest is to explore the structure of ACIPs
and the relationship between the uniqueness of such measures and
topological properties, e.g., the existence of dense orbits,
topological transitivity and minimality for multidimensional PAPs
(particularly for PWIs) by applying the functional analytic
approach. Definitions of PAPs and PWIs are given below.

Let $X$ be a compact subset of $\mathbb{R}^{d}$ and
$(X,\mathfrak{B},m)$ be a probability space. For convenience, $m$ always denotes $d-$dimension normalized Lebesgue measure on
$X$, and $\mathfrak{B}$ is the Borel $\sigma$-field. We say
$\mathcal{P}=\{\omega_{i}\}_{i=0}^{r-1}$ is a {\em topological
partition} of $X$ if: ({\bf i}) $\omega_{i}\cap
\omega_{j}=\emptyset,$ for $i\neq j;$~({\bf ii})
$\bigcup_{i=0}^{r-1}\omega_{i}=X$; and ({\bf iii}) for each
$\omega_{i}$, $ \operatorname{int}(\omega_{i})\neq\emptyset$ and
$m(\partial \omega_{i})=0$. Here each $\omega_{i}$ is called an
atom; $\operatorname{int}A$ and $\partial A$ are the interior and
boundary of $A$ respectively.

\begin{definition}\label{def_PAPs}
A nonsingular map $f:(X,\mathfrak{B},m)\rightarrow
(X,\mathfrak{B},m)$ with a topological partition
$\mathcal{P}=\{\omega_{i}\}_{i=0}^{r-1}$ is called a {\em piecewise
area preserving map} (PAP) if
$f|_{\operatorname{int}(\omega_{i})}\in C^{1}$ for each $\omega_{i}$
and $|\det Df(x)|\equiv 1$ for $x\in
\bigcup_{i=0}^{r-1}\operatorname{int}(\omega_{i})$. Here
non-singularity means that $f$ is measurable (with respect to
$\mathfrak{B}$) and $m(A)=0$ implies $m(f^{-1}(A))=m(f(A))=0$ for
any $A\in\mathfrak{B}$;  and $Df$ refers to the Jacobian matrix. We
say a PAP $f$ is piecewise-invertible-area-preserving if
$f|_{\omega_{i}}$ is invertible for each $\omega_{i}$, and say $f$
is an invertible PAP if $f$ is globally invertible. In particular,
if each $f|_{\operatorname{int}(\omega_{i})}$ is isometry (i.e.,
preserving Euclidean distance) then we say $f$ is a {\em piecewise
isometry} (PWI).
\end{definition}

Our definition of PAPs include piecewise hyperbolic maps with determinant~$\pm 1$,
e.g., baker's map, Arnold's cat map, area preserving H\'enon map and standard map \cite{Mackay93}.
However, we will mainly concentrate on methods working on non-hyperbolic maps such as PWIs.

For a PAP $f:X\rightarrow X$, the ACIPs are classified based on the
density properties as \footnote{$\frac{d\mu}{dm}\in L^{1}(m)$ is
$m-a.e.$ continuous means that its equivalence class contains an
$m-a.e.$ continuous representative.}
\begin{eqnarray*}
\mathcal{M}_{I}(f): &=& \{\mu~\mbox{is an ACIP with respect to $f$}\},\\
\mathcal{M}_{IB}(f): &=& \{\mu\in\mathcal{M}_{I}(f): \frac{d\mu}{dm}=\eta|_{X} \mbox{~for some~} \eta\in BV(\Omega), \mbox{where~}  \Omega\supset X \mbox{~is an~open ball}\},\\
\mathcal{M}_{IC}(f):&=&\{\mu\in\mathcal{M}_{I}(f):~\frac{d\mu}{dm}~\mbox{is~}m-a.e.\mbox{~continuous}\},
\end{eqnarray*}
where $BV(\Omega)$ is the space of bounded variation functions (see
Definition~\ref{def_bv}). We chose to work with $\mathcal{M}_{IB}$
and $\mathcal{M}_{IC}$ for the following reasons.
\begin{itemize}
  \item These spaces are ``large enough'' Banach subspaces of $L^{1}(m)$, i.e., they contain discontinuous functions \cite{Saussol00}.
  \item Functions in these spaces have ``good'' geometric properties, e.g., $\chi_{E}\in BV(\Omega)$ implies that the measurable
  subset $E\subset \Omega$ has finite perimeter \cite{Evans92}.
  \item These spaces coincide with those chosen in piecewise hyperbolic maps in
  \cite{Baladi09,Buzzi01,Tsujii00}.
  \item These spaces are invariant under the {\em transfer operator} (defined in Section~\ref{subsec_Perron}) for PWIs.
\end{itemize}

It is clear that $\mathcal{M}_{IB}\subset\mathcal{M}_{IC}\subset
\mathcal{M}_{I}$ for one dimensional invertible PAPs, while in
higher dimensions,
$\mathcal{M}_{IB}\cup\mathcal{M}_{IC}\subset\mathcal{M}_{I}$.
Additionally, for non-invertible PAPs, the set $\mathcal{M}_I$ is
possibly empty and conditions for which
$\mathcal{M}_{I}\neq\emptyset$ are discussed in Section
\ref{sec_pie_invertible}.

The novelty of this article is that we introduce multidimensional
bounded variation functions to analyze ACIPs for PAPs, especially for
PWIs. In Theorem~\ref{theor_1}, we explore the relationship between
the set of nomadic points and the sets $\mathcal{M}_{IC}$,
$\mathcal{M}_{IB}$ for invertible multidimensional PAPs. In
particular, we demonstrate that when the set of nomadic points has a
positive Lebesgue measure, both $\mathcal{M}_{IB}$ and $\mathcal{M}_{IC}$ are singletons. This can be applied to non-uniquely ergodic IETs constructed in \cite{Keane77,Keynes76} to show the irregularity of densities of their ACIPs. For invertible PWIs, in Theorem~\ref{theo_2} we give an approach to construct invariant measures with bounded variation densities. These results partially answer one of Goetz's questions in~\cite{Goetz03}.

The paper is organized in the following way. Preliminaries and the
main results are stated in Section~\ref{sec_def}, then applications
along with discussions are in Section~\ref{sec_app} and finally
proofs are given in Section \ref{sec_proof}.

\section{Preliminaries and Main results}\label{sec_def}
In this section, we give the formal definitions of {\em transfer
operator} and {\em multidimensional bounded variation}, and then
state the main results which are connected to one of the open
questions in \cite{Goetz03}.

\subsection{Transfer operator}\label{subsec_Perron}
Let $(X,\mathfrak{B},m)$ be a probability space where $m$ is
normalized Lebesgue measure and let $f:X\rightarrow X$ be a
nonsingular map. The {\em transfer operator}
$\mathcal{L}_{f}:L^{1}(m)\to L^{1}(m)$ associated with $f$ is
defined up to $m-a.e.$ equivalence as follows \cite{Boy97}:
\begin{equation*}\label{equ_def_Perron operator}
   \int_{A}\mathcal{L}_{f}\varphi dm=\int_{f^{-1}(A)}\varphi
   dm,~\varphi\in L^{1}(m),~A\in \mathfrak{B}.
\end{equation*}

This {\em transfer operator} possesses the following dual property \cite{Viana97}
\begin{equation*}\label{pro_dual}
\int(\mathcal{L}_{f}\varphi)\psi dm=\int\varphi\cdot(\psi\circ
f)dm,~~\varphi\in L^{1}(m),\psi\in L^{\infty}(m).
\end{equation*}

For an invertible PAP $f:X\rightarrow X$ with a topological
partition $\mathcal{P}:=\{\omega_{0},\cdots,\omega_{r-1}\}$, the
{\em transfer operator} can be simplified to be
$$
\mathcal{L}_{f}(\varphi)=\varphi\circ f^{-1}, ~\forall \varphi\in
L^1(m).
$$

\subsection{Multidimensional bounded variation}\label{sec_bv}
There are various definitions of multidimensional bounded variation functions, e.g., see Appendix~\ref{app_BV} and \cite{Evans92,Volper85}.
These definitions can be reduced to the usual notation of bounded variation in one dimension; see~Appendix~\ref{app_BV}. We state one of these as follows.
\begin{definition} \cite{Evans92}\label{def_bv}
Let $\Omega$ be an open set of $\mathbb{R}^{d}$. A function $\eta\in
L^{1}(\Omega)$ is a {\em bounded variation} function ($\eta\in
BV(\Omega)$) if
\begin{equation}\label{equ_variance}
\operatorname{var}(\eta):=\sup\left\{\int_{\Omega}\eta\cdot\operatorname{div}\overrightarrow{\phi}~
dm: \overrightarrow{\phi}\in
C_{c}^{1}(\Omega,\mathbb{R}^{d}),||\overrightarrow{\phi}||_{\infty}\leq 1\right\}<\infty.
\end{equation}
Here $\overrightarrow{\phi}=(\phi_{i})_{i=1}^{d}$,
$\operatorname{div}\overrightarrow{\phi}=\sum_{i=1}^{d}\frac{\partial\phi_{i}}{\partial
x_{i}},$ $||\overrightarrow{\phi}||_{\infty}:=\sup_x|\overrightarrow{\phi}(x)|$, and $C_{c}^{1}(\Omega,\mathbb{R}^{d})$ is the set of
$\overrightarrow{\phi}\in C^{1}(\Omega,\mathbb{R}^d)$ with compact support. We define a norm on $BV(\Omega)$ by
$||\eta||_{BV}:=||\eta||_{1}+\operatorname{var}(\eta)$.
\end{definition}

For functions of bounded variation, we state the corresponding
Helly's Theorem \cite{Evans92} below.

{\bf Helly's Theorem} \cite{Evans92} {\em Let $\Omega\subset
\mathbb{R}^{d}$ be an open and bounded domain with Lipschitz boundary. Assume that $\{\eta_{n}\}_{n=1}^{\infty}$ is a
sequence in $BV(\Omega)$ satisfying $\sup_{n}||\eta_{n}||_{BV}<\infty$, then there exist a
subsequence $\{\eta_{n_{k}}\}_{k=1}^{\infty}$ and a function
$\eta\in BV(\Omega)$ such that $\eta_{n_k}\to \eta$ in
$L^{1}(\Omega)$ as $k\rightarrow\infty$.}

\subsection{Main results}\label{sec_mainresult}
Let $X$ be a compact subset of $\mathbb{R}^{d}$ and $f:X\rightarrow
X$ be an invertible map. We say $x\in X$ is a {\em nomadic point} of $f$ if $O_{f}(x):=\{f^{i}(x)|i\in
\mathbb{Z}\}$ is dense in $X$. We denote
$\operatorname{nom}(f)$ to be the set of all nomadic points of the
map $f$. If $\operatorname{nom}(f)=X$ then $f$ is called {\em
minimal}.
\begin{theorem}\label{theor_1}
Let $(X,\mathfrak{B},m)$ be a probability space where $m$ is the normalized Lebesgue measure and
$f: X\rightarrow X$ be an invertible PAP with a topological
partition $\mathcal{P}=\{\omega_{i}\}_{i=0}^{r-1}$. Then the following
hold:
\begin{enumerate}
\item[$(i)$] if $m(\operatorname{nom}(f))>0$ then $\mathcal{M}_{IB}(f)\cup \mathcal{M}_{IC}(f)=\{m\}$;

\item[$(ii)$] if $f|_{\omega_i}$ is a homeomorphism for each $\omega_i$, then $\mathcal{M}_{IC}(f)=\{m\}$
implies $\operatorname{nom}(f)\neq\emptyset$.
\end{enumerate}
\end{theorem}

We remark here that even for area preserving diffeomorphisms, $nom(f)\neq \emptyset$ does not necessarily imply $m(nom(f))>0$, see e.g. Fayad and Katok's \cite{Fayad04}.

\begin{corollary}\label{cor_1}
Suppose $f$ is an invertible PAP with $m(\operatorname{nom}(f))>0$ and there exists a measure $m\neq \mu\in\mathcal{M}_{I}(f)$, then
$\varphi:=\frac{d\mu}{dm}\notin BV(\Omega)$ and the set of
discontinuities of $\varphi$ has a positive Lebesgue measure.
\end{corollary}

The following Theorem~\ref{theo_2} aims to construct ACIPs with bounded variation densities for invertible PWIs, say $f:X\rightarrow X$ with a topological partition $\mathcal{P}=\{\omega_0,\cdots,\omega_{r-1}\}$.
As a bounded variation function is defined on an open set, we choose an open ball $\Omega\supset X$ and extend $f$ to
$\overline{f}:\Omega\rightarrow \Omega$ by
\begin{equation}\label{equ_fbar}
\overline{f}(x)=\left\{\begin{array}{ll}
          f(x),&x\in X\\
           x,& x\in \Omega\backslash X.
         \end{array}
          \right.
\end{equation}
Given any $\eta\in BV(\Omega)$, the sequence of variations
$\{\operatorname{var}(\mathcal{L}^n_{\overline{f}}\eta)\}_{n=0}^{\infty}$
are not necessarily uniformly bounded \cite{Keller96}. Therefore, we
alternatively work on functions $\eta$ which lie in a plausible
proper subset $BV^*(\Omega)$ (see below). This subset is associated
with a Sobolev space.

Let $\omega_r:=\Omega\backslash X$ and without ambiguity we still write $\mathcal{P}=\{\omega_0,\cdots,\omega_r\}$ as a topological partition of $\Omega$. Moreover, we denote
$$\partial\mathcal{P}^{\infty}:=\{x\in\Omega: \overline{f}^n(x)\in\partial\mathcal{P}\mbox{~for some~}n\geq 0\},
$$
where $\partial\mathcal{P}:=\cup_{i=0}^r\partial\omega_{i}$, and define a $\delta-$neighborhood of
$\partial\mathcal{P}^{\infty}$ by
\begin{equation}\label{equ_Nd}
N_{\delta}:=\{x\in\Omega,\operatorname{dist}(x,\partial\mathcal{P}^{\infty})<\delta\}.
\end{equation}
For a given invertible PWI $f$, the function subspace $BV^*(\Omega)$ that we consider is defined as
\begin{equation}\label{def_bvstar}
BV^*(\Omega):=\{\eta\in BV(\Omega):\eta|_{N_{\delta}}\in W^{1,2} \mbox{~for some~}\delta>0 \},
\end{equation}
where $W^{1,2}$ is a Sobolev space (see Appendix~\ref{app_sobo}). We
remark that different invertible PWI $f$ determines different
$BV^*(\Omega)$ individually, but in all cases
$W^{1,2}(\Omega)\subseteq BV^*(\Omega)\subseteq BV(\Omega)$.

\begin{theorem}\label{theo_2}
Suppose $f:X\rightarrow X$ is an invertible PWI and $\Omega\supset
X$ an open ball. Then given any $\eta\in BV^{*}(\Omega)$ with $\eta|_X\geq 0$ and $||\eta|_X||_1>0$, there exists a subsequence of the Birkhoff average of the transfer
operator $\mathcal{L}_{\overline{f}}$, which converges to a function
$\overline{\eta}\in BV(\Omega)$ in $L^1(m)$, i.e.,
$$\frac{1}{n_k}\sum_{i=0}^{n_k-1}\mathcal{L}_{\overline{f}}^{i}\eta\rightarrow \overline{\eta}\in BV(\Omega), \mbox{~as~}k\rightarrow \infty,$$
and by normalization $d\mu:=\overline{\eta}|_{X}dm\in \mathcal{M}_{IB}(f)$.
\end{theorem}


Concerning the open question \cite{Goetz03} in piecewise rotations (defined in Appendix~\ref{app_piecerotation}), the following corollaries give a universal approach to partially determine the ACIPs.

\begin{corollary}\label{col_rotation_inv}
Suppose $f: X\rightarrow X$ is an invertible piecewise rotation, then
\begin{enumerate}
  \item[$(i)$] $\mathcal{M}_{I}(f)=\{\varphi dm: \varphi=\mathbb{E}(\varphi|\mathcal{I}), \varphi\in L^1(m)\}$, where $\mathcal{I}=\{B\in\mathcal{B}: f^{-1}(B)=B \mod m\}$;
  \item[$(ii)$] given any $\eta\in BV^*(\Omega)$ satisfying $\eta|_X\geq 0$ and $||\eta|_X||_1>0$, any accumulation point of $\{\frac{1}{n}\sum_{i=0}^{n-1}\mathcal{L}_{\overline{f}}^i\eta\}_{n=1}^{\infty}$ is an invariant density of the map $f$. Furthermore, if $m(\operatorname{nom}(f))>0$ then $\mathcal{M}_{IB}(f)=\{m\}$.
\end{enumerate}
\end{corollary}

\begin{corollary}\label{col_rotation_noninv}
Suppose $f:X\rightarrow X$ is a non-invertible piecewise rotation. Let $X^+:=\bigcap_{i=0}^{\infty} f^i(X)$ and define $f^+:\overline{X^+}\rightarrow\overline{X^+}$ as in equation~\eqref{equ_fplus} in Section~\ref{sec_pie_invertible}. Then $f^+$ is $m-a.e.$ invertible. Furthermore,
\begin{itemize}
 \item[$(i)$] if $m(\overline{X^+})>0$, then the statements in Corollary~\ref{col_rotation_inv} hold for $f^+$ and
$$
\mathcal{M}_{I}(f)=\{\mu(\cdot):=\nu(\cdot\cap\overline{X^+}),\forall\nu\in\mathcal{M}_{I}(f^+)\};
$$
\item[$(ii)$] if $m(\overline{X^+})=0$, then $\mathcal{M}_{I}(f)=\mathcal{M}_{IB}(f)=\mathcal{M}_{IC}(f)=\emptyset.$
\end{itemize}
\end{corollary}

Proof of Corollary~\ref{col_rotation_inv} is based on Theorem~\ref{theor_1}, \ref{theo_2} and Lemma~\ref{lem_Condition_expectation} while proof of Corollary~\ref{col_rotation_noninv} is based on Lemma~\ref{lem_almostclosedness} and Proposition~\ref{pro_nonin}. We remark that ACIPs only give a subset of non-atomic probability Borel measures. Therefore, to fully answer the question in piecewise rotations \cite{Goetz03}, we have to explore singular non-atomic probability invariant measures. For instance, when $m(\overline{X^+})=0$, it is natural to consider Hausdorff measure. This is discussed at the end of Section~\ref{sec_pie_invertible}.

\section{Applications and Discussions}\label{sec_app}
In this section, we consider two main applications. We first consider IETs (see Appendix~\ref{app_IET} for the definition) which are non-uniquely ergodic. We apply Theorem~\ref{theor_1} to show that the densities of their ACIPs can be irregular. We then consider multidimensional piecewise invertible area preserving maps and apply Theorem~\ref{theor_1}
and Theorem~\ref{theo_2} to study their invariant densities. At the end of this section, we give a short discussion on the open question posed in \cite{Goetz03} for piecewise rotations.

\subsection{Interval exchange transformations}\label{sec_IETs}
For an IET $f$, the set $\mathcal{M}_{IC}(f)$ can be refined to
\[
\mathcal{M}_{IC}'(f):=\{\mu\in\mathcal{M}_{I}(f):~\frac{d\mu}{dm}:=\varphi~\mbox{has~at~most~countably
many~discontinuous~points}\},
\]
where $m$ is the normalized Lebesgue measure. Observe that
$\mathcal{M}_{IB}(f)\subset\mathcal{M}^{'}_{IC}(f)\subset\mathcal{M}_{IC}(f)$. Moreover, it is known that topological transitivity\footnote{Topological transitivity
means that for any open sets $U$ and $V$, there exists $n\in
\mathbb{Z}$ such that $f^{n}(U))\cap V\neq\emptyset$.} implies
minimality for IETs (see e.g., Corollary $14.5.11$ in
\cite{Katok95}). Hence by applying Theorem \ref{theor_1}, we have the following corollary which characterizes
the minimality properties of IETs.
\begin{corollary}\label{coro_IET}
For any IET $f:[0,1)\rightarrow [0,1)$,
\[
f\mbox{~is
minimal}\Leftrightarrow\mathcal{M}_{IC}(f)=\{m\}\Leftrightarrow\mathcal{M}'_{IC}(f)=\{m\}.
\]
\end{corollary}

This corollary can be used to investigate Keane's conjecture, namely that minimality implies unique ergodicity for IETs \cite{Keane75}. This conjecture was shown to be false and we review two well-known counterexamples here.

Keynes and Newton \cite{Keynes76} considered the following map. Let
$\widehat{T}_{\gamma\beta}:[0,1+\beta)\rightarrow[0,1+\beta)$ to be
\[
\widehat{T}_{\gamma\beta}(x)=\left\{\begin{array}{ll}
                     x+1, & \mbox{if}~ 0\leq x<\beta \\
                     x+\gamma~(\mbox{mod}~1), & \mbox{if}~\beta\leq
                     x<1+\beta.
                   \end{array}\right.
\]
By choosing appropriate $\beta$ and $\gamma$ (see
\cite{Coffey88,Keynes76}), the map
$T_{\gamma\beta}(x):=\frac{1}{1+\beta}\widehat{T}_{\gamma\beta}(x(1+\beta))$
is minimal and has an eigenvalue $-1$. This implies that $T_{\gamma\beta}^{2}$ is not
uniquely ergodic and its ergodic measures belong to $\mathcal{M}_{I}$
(see \cite{Keynes76} for details).

Keane \cite{Keane77} also constructed an IET with four intervals satisfying a strong irrationality condition that implies minimality. Under certain conditions, there exist two different ergodic
measures $\mu_{1}$ and $\mu_{2}$. Moreover, such ergodic measures
are either both in $\mathcal{M}_{I}$ or one is Lebesgue measure and
the other is singular. For the measure that is singular, the Hausdorff dimension has been recently estimated \cite{Chaika08}. Together with our results, we could obtain a better understanding of ergodic measures for non-uniquely ergodic IETs.

For the examples above, there are no explicit formulae for their
densities (even if these densities belong to $\mathcal{M}_{I}$). One
of the difficulties in constructing counterexamples can be seen from
the fact that these ergodic measures are in
$\mathcal{M}_I(f)\backslash\mathcal{M}_{IC}(f)$ from
Corollary~\ref{coro_IET}. In the following proposition, we provide
a more explicit description of the invariant densities for non-uniquely
ergodic IETs.

\begin{proposition}\label{prop_density}
Let $f$ be any topologically transitive IET on $[0,1)$. Suppose $m\neq
\mu\in\mathcal{M}_{I}(f)$, then the following hold:
\begin{enumerate}
  \item[$(i)$] the density of $\mu$ is a simple function (i.e., a linear combination of finitely many characteristic functions);
  \item[$(ii)$] for any representative from the equivalence class
  $\varphi:=\frac{d\mu}{dm}$, $\varphi$ is discontinuous everywhere and $\operatorname{supp}\mu=[0,1)$
  (recall that if $x\in\operatorname{supp}\mu$, then for any open ball $B_x$ containing $x$, $\mu(B_x)>0$ hold \cite{Katok95}).
\end{enumerate}
\end{proposition}
\begin{remark}
For the two ergodic measures $\mu_{1},\mu_{2}\in\mathcal{M}_{I}$
(i.e., $\mu_1,\mu_2\ll m$) in the examples of Keynes $\&$
Newton~\cite {Keynes76} and Keane \cite{Keane77}, we can derive that
$\mu_{1}\perp\mu_{2}$ \cite{Katok95} and moreover,
$\operatorname{supp}\mu_{1}=\operatorname{supp}\mu_{2}=[0,1)$ from
Proposition \ref{prop_density}. Hence, the measures $\mu_{1}$ and
$\mu_{2}$ intermingle with each other in some sense.
\end{remark}

\subsection{Piecewise invertible area preserving maps}\label{sec_pie_invertible}

In this subsection, we aim to understand the structure of ACIPs for piecewise invertible area preserving maps $f:X\rightarrow X$. For such a map $f$, the set $X^+:=\bigcap_{i=0}^{\infty}f^{i}(X)$ is invariant under the map, i.e., $f(X^+)=X^+$ \cite{Fu07pa}. In particular for PWIs, $X^{+}$ is almost closed, i.e.,
$m(X^{+})=m(\overline{X^{+}})$ \cite{Ashwin02}. Here we show that such almost closedness of $X^+$ is valid for a broad range of piecewise invertible area preserving maps.

\begin{lemma}\label{lem_almostclosedness}
\label{lemma_piece_invar} Let $f:X\rightarrow X$ be a piecewise
invertible area preserving map with a topological partition
$\mathcal{P}=\{\omega_{0},\cdots,\omega_{r-1}\}$. Suppose
$f_i:=f|_{\operatorname{int}\omega_i}$ is Lipschitz for each
$\omega_i$, then $m(X^+)=m(\overline{X^+})$ and $f|_{X^+}$ is
$m-a.e.$ invertible.
\end{lemma}

Under the conditions of Lemma~\ref{lem_almostclosedness}, it is not necessary to have the property
$f(\overline{X^{+}})\subseteq \overline{X^+}$, but we can define a map $f^{+}$ that is $m-a.e.$ equal to $f$ and for which $f^+(\overline{X^{+}})\subset\overline{X^+}$  (see below). Since each $f_i$ is Lipschitz, then there exists a continuous extension $\widehat{f_{i}}:\overline{\operatorname{int}\omega_{i}}\to \overline{f_{i}(\operatorname{int}\omega_{i})}$. For any $x\in(\bigcup_{i=0}^{r-1}\partial\omega_{i})\cap\overline{X^{+}}$,
let $g(x):=\widehat{f_{i^*}}(x)$  where
$i^{*}:=\min\{i:x\in\partial\omega_{i}\}$. Then we can define
$f^+:\overline{X^+}\rightarrow \overline{X^+}$ to be
\begin{equation}\label{equ_fplus}
f^+(x)=\left\{
\begin{array}{lr}
f(x), ~x\in\operatorname{int}(\omega_i)\cap X^+\\
g(x), ~\mbox{otherwise}.
\end{array}
\right.
\end{equation}
Moreover, if $f_i$ is bi-Lipschitz, the map $f^+$ can be shown to be
non-singular and to be $m-a.e.$ invertible. The non-singularity and
$m-a.e.$ invertibility of $f^+$ allow to obtain
Theorem~\ref{theor_1} and Theorem~\ref{theo_2} for
$f^+:\overline{X^+}\rightarrow \overline{X^+}$. 
The ACIPs of $f^+$ can further be used to determine the ACIPs of $f$ as follows.
\begin{proposition}\label{pro_nonin}
Let $f:X\rightarrow X$ be a piecewise invertible area preserving map
with a topological partition
$\mathcal{P}:=\{\omega_{0},\cdots,\omega_{r-1}\}$ and $f^{+}:\overline{X^{+}}\to\overline{X^{+}}$ be defined
as in equation~\eqref{equ_fplus}. Suppose that
$f_i:=f|_{\operatorname{int}\omega_i}$ is bi-Lipschitz continuous. Then $f^+$ is non-singular and $m-a.e.$
invertible. Moreover, the following hold:
\begin{enumerate}
  \item [$(i)$] if $m(\overline{X^{+}})>0,$ then $\mathcal{M}_{I}(f)=\left\{\mu(\cdot):=\nu(\cdot\cap\overline{X^+}),\forall\nu\in\mathcal{M}_{I}(f^+)\right\}$;
  \item [$(ii)$] if $m(\overline{X^{+}})=0,$ then $\mathcal{M}_{I}(f)=\emptyset.$
\end{enumerate}
\end{proposition}

When $m(\overline{X^{+}})=0$, it is natural to consider invariant
measures that are absolutely continuous with respect to Hausdorff
measure. If we let $s=\dim_{H}\overline{X^{+}}$ and furthermore, if $f^+$ satisfies the following conditions:
\begin{description}
\item [(1)] $0<\mathcal{H}^{s}(X^{+})=\mathcal{H}^{s}(\overline{X^{+}})<\infty;$
\item [(2)] $\mathcal{H}^{s}$ is an invariant measure for $f^{+};$
\item [(3)] $f^+$ is non-singular with respect to $\mathcal{H}^s$, i.e.,
$\mathcal{H}^{s}\left((f^{+})^{-1}(A)\right)=\mathcal{H}^{s}(f^+(A))=0$
whenever $\mathcal{H}^{s}(A)=0$;
\end{description}
then by the arguments analogous to those used in Proposition \ref{pro_nonin}, we can show that
$f^{+}$ is $\mathcal{H}^{s}-a.e.$ invertible.

The above three conditions can be achieved for some piecewise invertible area preserving maps. We take interval translation maps (see Appendix~\ref{app_IET}) as examples. Condition (2) is demonstrated in~\cite{Bruin03} while condition (3) can be inferred by combining condition (2) and the definition of Hausdorff measure. The $\mathcal{H}^{s}-a.e.$ closedness of $X^+$ can be shown by adapting the proof of Lemma~\ref{lemma_piece_invar}. Moreover, by \cite[Theorem 9.3]{Falconer03}, condition~(1) hold for particular interval translation maps where $\overline{X^+}$ are self similar sets satisfying an open set condition and positive Hausdorff dimension \cite{Bruin03}.

Concerning the open question in \cite{Goetz03}, for non-invertible piecewise rotations in the case of $m(\overline{X^+})=0$, we consider absolutely continuous (with respect to $\mathcal{H}^s$) invariant probability measures. We denote
$$
\mathcal{H}^{s}_{I}(f):=\{\nu~\mbox{probability invariant measure of
f}~:\nu\ll\mathcal{H}^{s}\}.
$$
\begin{proposition}\label{prop_hausdorff}
Suppose $f:X\rightarrow X$ is a two-dimensional piecewise rotation with
$m(\overline{X^{+}})=0$ and $s:=\dim_{H}\overline{X^{+}}>1$, then
$X^{+}$ is $\mathcal{H}^{s}-a.e.$ closed and $f^+$ is non-singular with respect to $\mathcal{H}^s$.
Moreover, the following hold:
\begin{enumerate}
  \item[(i)] if $0<\mathcal{H}^{s}(\overline{X^{+}})<\infty$, then $\mathcal{H}^{s}$ is an invariant measure of $f^{+}$ and $f^{+}$ is $\mathcal{H}^{s}-a.e.$ invertible; furthermore, $\mathcal{H}^s_I(f)=\left\{\mu(\cdot):=\nu(\cdot\cap \overline{X^+}),\forall\nu\in\mathcal{H}^s_{I}(f^{+})\right\}$;
  \item[(ii)] if $\mathcal{H}^{s}(\overline{X^{+}})=0$ then $\mathcal{H}^s_{I}(f)=\emptyset$.
\end{enumerate}
\end{proposition}

Under the condition $s:=\dim_{H}\overline{X^{+}}>1$, the proof of Proposition~\ref{prop_hausdorff} is analogous to the proofs of Proposition~\ref{pro_nonin} and Lemma~\ref{lemma_piece_invar}. However, this condition does not always hold. For instance, the Cartesian product of interval translation maps (as defined in \cite{Bruin03}) with themselves provides some examples of piecewise rotations with Hausdorff dimension ranging $0\leq s\leq 1$. It might be interesting to explore conditions for $s>1$.

Regarding the determination of all the invariant non-atomic probability measure for piecewise rotations, it might be necessary to consider the structure of $\mathcal{H}^s_I(f^+)$ in the case of $\mathcal{H}^{s}(\overline{X^{+}})=\infty$. In this case, $f^+$ is not necessarily $\mathcal{H}^s-a.e.$ invertible, however, we note that by \cite[Theorem 6.2]{Falconer03}, there exists a compact subset $E\subset \overline{X^{+}}$ with $0<\mathcal{H}^{s}(E)<\infty.$ We suggest to establish a non-atomic probability invariant measures of $f^{+}$ induced by $E$ as a reference measure and leave this for further studies.

\section{Proofs}\label{sec_proof}
We first state a basic lemma regarding $\mathcal{M}_{I}(f)$, i.e., the set of all ACIPs with respect to a map $f$, with a standard proof.

\begin{lemma}\label{lemma_fixed}
Let $f:X\to X$ be a nonsingular map and $\varphi\in L^{1}(m)$, then
$\mathcal{L}_f\varphi=\varphi$ if and only if $d\mu:=\varphi dm\in
\mathcal{M}_{I}(f)$.
\end{lemma}
{\bf \  \ Proof:} Suppose $\mathcal{L}_f\varphi=\varphi$. Let
$d\mu:=\varphi dm$, then for each $A\subset X$,
 \[
\mu(A)=\int_{A}\varphi dm=\int_{A}\mathcal{L}_{f}\varphi
dm=\int_{f^{-1}(A)}\varphi dm=\mu(f^{-1}(A)),
 \]
which implies $\mu\in \mathcal{M}_{I}(f)$.

On the other hand, if $d\mu:=\varphi dm$ is an invariant measure of $f$, i.e., $\mu(f^{-1}(A))=\mu(A)$ for any Borel set $A\subset X$, then $\int_{A}\mathcal{L}_{f}\varphi dm=\int_{f^{-1}(A)}\varphi dm=\mu(f^{-1}(A))=\mu(A)=\int_{A}\varphi dm $, which implies $\mathcal{L}_{f}\varphi=\varphi$.  \hfill $\Box$

\subsection{Proof of Theorem \ref{theor_1}}
\label{proof_1}
To prove statement $(i)$, we first show that $\mathcal{M}_{IC}=\{m\}$ followed by a proof of $\mathcal{M}_{IB}=\{m\}$. 
For statement $(ii)$, we start with a lemma showing the equivalence between topological transitivity and the existence of a nomadic point.

{\bf \ \ Proof of the statement (i) in Theorem \ref{theor_1}:} Consider any $\mu\in \mathcal{M}_{IC}$. Take a representative $\varphi=\frac{d\mu}{dm}$ from the equivalence class such that $\varphi$ is $m-a.e.$ continuous. Furthermore, take any point $x'\in\bigcup_{i=0}^{r-1}\operatorname{int}\omega_i$ where $\varphi$ is continuous. Since $m(nom(f))>0$, we can choose a nomadic point $x^*$ such that for any $n\in \mathbb{Z}$, the equality $\varphi\circ f^{-n}(x^*)=\varphi(x^*)$ holds. Then there exists a subsequence
$\{f^{k_t}(x^{*})\}$ such that as $|k_t|\to\infty$, $f^{k_t}(x^{*})\to
x'$ because $x^*$ is a nomadic point. By the continuity of $\varphi$ at $x'$, we have $$\varphi(x')=\varphi\left(\lim_{|k_t|\to \infty}f^{k_t}(x^*)\right)=\varphi(x^*).$$ This implies $\varphi\equiv 1$, i.e., $\mathcal{M}_{IC}=\{m\}$.

Consider the case $\mu\in\mathcal{M}_{IB}(f)$, i.e.,
$\varphi=\frac{d\mu}{dm}$ with $\varphi=\eta|_{X}$ for some $\eta\in
BV(\Omega)$. Hence, by \cite[page 178]{Volper85}, $m-a.e.~x\in\Omega$ are {\em regular points} of $\eta$. By a regular point $x_0$, we mean there exists a unit vector $a\in\mathbb{R}^d$ such that the limits
$$\lim\limits_{x\rightarrow x_0,\left<x-x_0,a\right>>0}\eta(x)\mbox{~~and~}\lim\limits_{x\rightarrow x_0,\left<x-x_0,a\right><0} \eta(x)
$$ exist, where $\left<\cdot,\cdot\right>$ is the inner product. Therefore, by analogous arguments to that used in $\mathcal{M}_{IC}(f)=\{m\}$, it follows that for any regular point $x_0\in\bigcup^{r-1}_{i=0}\operatorname{int}(\omega_i)$,
$$
\lim\limits_{x\rightarrow x_0,\left<x-x_0,a\right>>0}\varphi(x)=\lim\limits_{x\rightarrow x_0,\left<x-x_0,a\right><0} \varphi(x).
$$
Hence, by \cite[page 168]{Volper85}, we have $\lim_{x\to x_{0}}\varphi(x)=\varphi(x^{*})$ where $x^*$ is a nomadic point. In addition, $\varphi\in L^{1}(m)$ implies that $m-a.e.$ $x\in X$ is a Lebesgue point of $\varphi$, i.e.,
\[
\lim_{r\to
0^{+}}\frac{1}{m(B(x,r))}\int_{B(x,r)}|\varphi(y)-\varphi(x)|dm(y)=0.
\]
Therefore, if a regular point $x_0$ is also a Lebesgue point, then
\begin{eqnarray*}
0\leq|\varphi(x_0)-\varphi(x^{*})|=\lim_{r\to
0^{+}}\frac{1}{m(B(x_0,r))}\int_{B(x_0,r)}|\varphi(x_0)-\varphi(x^{*})|dm(y)\\
\leq\lim_{r\to
0^{+}}\frac{1}{m(B(x_0,r))}\int_{B(x_0,r)}|\varphi(y)-\varphi(x_0)|dm(y)\\
+\lim_{r\to
0^{+}}\frac{1}{m(B(x_0,r))}\int_{B(x_0,r)}|\varphi(y)-\varphi(x^{*})|dm(y)=0.
\end{eqnarray*}
This implies that $\varphi(x_0)=\varphi(x^{*})$, meaning that $\varphi\equiv
1$. \hfill$\Box$

Before proving statement (ii), we first formulate an equivalent condition of topological transitivity.

\begin{lemma}\label{lem_transitive}
Let $f:X\rightarrow X$ be an invertible PAP with a
topological partition $\mathcal{P}=\{\omega_{i}\}_{i=0}^{r-1}$.
Suppose $f|_{\omega_{i}}$ is a homeomorphism for each $\omega_{i}$,
then the following are equivalent.
\begin{enumerate}
  \item [$(i)$] $f$ has a nomadic point;
  \item [$(ii)$] $f$ is strongly topologically transitive, i.e., for any open sets $U,V$, there exists $n\in \mathbb{Z}$ such that
$\mbox{int}(f^{n}(U))\cap V\neq\emptyset$;
  \item[$(iii)$] $f$ is topologically transitive.
\end{enumerate}
\end{lemma}

Since PAPs are not necessarily continuous at every point,
the proof of Lemma~\ref{lem_transitive} will not be standard (see \cite{Walters82} for the continuous version). Therefore, we provide the details of the proof here.

{\bf \ \ Proof of Lemma~\ref{lem_transitive}:} ``$(ii)$ implies $(iii)$'' is direct and we only need to prove $(i)$ implies $(ii)$ and $(iii)$ implies $(i)$. For convenience, we denote $\mathcal{P}^{(n)}:=\bigvee_{i=0}^{n}f^{-i}(\mathcal{P})$ for $n\geq 0$ and $\mathcal{P}^{(n)}:=\bigvee_{i=-1}^{n}f^{-i}(\mathcal{P})$ for $n<0$, and denote $\omega^{(n)}$ one of the atoms in the topological partition
$\mathcal{P}^{(n)}$ for any $n\in \mathbb{Z}$.

``$(i)\Rightarrow(ii)$.'' We prove by contradiction.
Suppose there exist open sets $U,V\neq\emptyset$ such that for any $n\in \mathbb{Z}$, $\operatorname{int}(f^{n}(U))\cap V=\emptyset$. Given a nomadic point $x^{*}$ of $f$, there exist $n_{1},n_{2}\in \mathbb{Z}$
 such that $f^{n_{1}}(x^{*})\in U$ and $f^{n_{2}}(x^{*})\in V$. Let $t:=n_{2}-n_{1}$, then there exists an atom $\omega^{(t)}\in\mathcal{P}^{(t)}$ such that $f^{n_1}(x^*)\in\omega^{(t)}$.

Suppose $f^{n_{1}}(x^{*})\in \operatorname{int}\omega^{(t)}$. Since $f^{t}|_{\operatorname{int}\omega^{(t)}}$ is homeomorphic, we have
$f^{t}(U\cap\operatorname{int}\omega^{(t)})=\operatorname{int}f^{t}(U\cap\operatorname{int}\omega^{(t)})$. Therefore,
\[
f^{n_{2}}(x^{*})=f^{t}(f^{n_{1}})(x^{*})\in
f^{t}(U\cap\operatorname{int}\omega^{(t)})\subset\operatorname{int}f^{t}(U),
\]
which implies $f^{n_{2}}(x^{*})\in \operatorname{int}(f^{t}(U))\cap
V.$ This is a contradiction.

Suppose $f^{n_1}(x^*)\in\partial\omega^{(t)}\cap\omega^{(t)}$. Let $l:=n_{1}+n_{2}$,
then there exists an atom $\omega^{(l)}$ such that $x^*\in
\omega^{(l)}$. Since
$f^{n_{1}}|_{\omega^{(l)}},f^{n_{2}}|_{\omega^{(l)}}$ are
continuous, there exists $x'\in\omega^{(l)}$ sufficiently
close to $x^{*}$ such that
$f^{n_{1}}(x')\in\operatorname{int}\omega^{(t)}\cap U$ and
$f^{n_{2}}(x')\in V.$ Repeat the same process as the above case,
using $x'$ in place of $x^{*}$ and this completes the proof.

``$(iii)\Rightarrow (i)$.'' Suppose $\{U_{i}\}_{i=1}^{\infty}$ is a
countable base for $X$. By $(iii)$, there exists an $n_{1}\in \mathbb{Z}$ such that
$f^{n_{1}}(U_{1})\cap U_{2}\neq\emptyset$. Let $y:=f^{n_{1}}(x)\in f^{n_{1}}(U_{1})\cap U_{2}$ where $x\in\omega^{(n_{1})}\in\mathcal{P}^{(n_{1})}$.

If $x\in\operatorname{int}\omega^{(n_{1})}$, then there exists an open ball $x\in B_{x}\subset\omega^{(n_{1})}$ such that $f^{n_{1}}(B_{x})\subset \operatorname{int}f^{n_1}(U_{1})$. Therefore, $y\in\operatorname{int}(f^{n_{1}}(U_{1}))\cap U_{2}$.

Otherwise, if $x\in\partial\omega^{(n_{1})}\cap\omega^{(n_1)}$, by using the approach analogous
to that used in ``$(i) \Rightarrow (ii)$'', then there exists $x'$ such that $x'\in\operatorname{int}\omega^{(n_{1})}\cap U_1 $ and
$y':=f^{n_{1}}(x')\in\operatorname{int}(f^{n_{1}}(U_{1}))\cap U_{2}.$ Hence
$$\operatorname{int}(f^{n_{1}}(U_{1}))\cap U_{2}\neq\emptyset.$$
Therefore, there exists a closed ball $B_{2}$ such that
$B_{2}\subset\operatorname{int}(f^{n_{1}}(U_{1}))\cap
U_{2}\cap\operatorname{int}\omega^{(n_{1})}.$ Moreover, since
$f^{n_1}|_{B_2}$ is a homeomorphism, then $V_{1}:=f^{-n_{1}}(B_2)$ is closed. Analogously, for
open sets $\operatorname{int}B_2$ and $U_{3}$ there exist
$n_{2}\in\mathbb{Z}$ and a closed ball $B_{3}\subset
\operatorname{int}(f^{n_2}(B_2))\cap U_3$, with $f^{n_{2}}|_{B_{3}}$ being
a homeomorphism. Let $V_{2}:=f^{-n_{2}}(B_3)\subset B_2$, then
$V_{2}$ is closed and $f^{-n_1}(V_2)\subset V_1$.

If one continues this process, there will exist
$\{n_{i}\}_{i=1}^{\infty}$ and a sequence of nonempty closed sets
$\{V_{i}\}_{i=1}^{\infty}$ such that $f^{-n_{i}}(V_{i+1})\subset
V_{i}$ for each $i$. Therefore,
$\bigcap_{i=1}^{\infty}f^{N}(V_{i+1})\neq\emptyset,$ where
$N=-\sum_{j=1}^{i}n_{j}.$ We note that $V_i\subset U_i$ for each $i$. Fix any
$\bar{x}\in\bigcap_{i=1}^{\infty}f^{N}(V_{i+1})$, then $\bar{x}$ is
nomadic. This completes the proof. \hfill $\Box$
\begin{remark}
There does exist a map $f$ that is topologically transitive but has no nomadic points~\cite{Peris99}. However, Lemma~\ref{lem_transitive} does not apply to this case as $f$ does not extend continuously to the boundary from its interior.
\end{remark}

{\bf Proof of statement (ii) in Theorem \ref{theor_1}:} We show by contradiction. Suppose $f$ has no nomadic points, so by Lemma \ref{lem_transitive} there will exist two open sets $U,V\subset X$ such that $\operatorname{int}(f^{n}(U))\cap V=\emptyset$ for all $n\in \mathbb{Z}$. Therefore,
$\operatorname{int}(f^{i+n}(U))\cap
\operatorname{int}f^{i}(V)=\emptyset~\operatorname{mod}~m,$ for any
$n,i\in \mathbb{Z}.$ Let
\[
U^{*}:=\bigcup_{i=-\infty}^{\infty}\operatorname{int}f^{i}(U)\mbox{~and~}V^{*}:=\bigcup_{i=-\infty}^{\infty}\operatorname{int}f^{i}(V).
\]
Then $U^{*}\cap V^{*}=\emptyset~\operatorname{mod}~m$ and both
$U^{*}$ and $V^{*}$ are invariant under $f$ up to $m-a.e.$. Hence both $m|_{U^{*}}$
and $m|_{V^{*}}$ are invariant measures and $m|_{U^{*}}\neq
m|_{V^{*}}.$ Since $\partial
U^{*}\subset\bigcup_{i=-\infty}^{\infty}\partial f^{i}(U)$ and
$\partial V^{*}\subset\bigcup_{i=-\infty}^{\infty}\partial
f^{i}(V)$, then this implies $m(\partial U^{*})=m(\partial
V^{*})=0$. I.e., the discontinuous points of $\chi_{U^{*}}$ and
$\chi_{V^{*}}$ lie in a set of zero Lebesgue measure. Therefore,
$m|_{U^{*}},m|_{V^{*}}\in \mathcal{M}_{IC}$, which contradicts the
uniqueness of measures in $\mathcal{M}_{IC}.$\hfill$\Box$

\begin{remark}
When restricting PAPs to IETs, each open set of $X$ is a union of countably many open intervals, therefore, $\chi_{U^{*}}$ and $\chi_{V^{*}}$ have at
most countably many discontinuities. This fact is used to prove Corollary~\ref{coro_IET}.
\end{remark}

\subsection{Proof of Theorem~\ref{theo_2}}
We prove Theorem~\ref{theo_2} by the following consecutive lemmas.
Recall that given an invertible PWI $f: X\rightarrow X$ and an open
ball $\Omega\supset X$, we can extend $f$ into an
$\overline{f}:\Omega\to\Omega$ as in \eqref{equ_fbar} and define
$BV^*(\Omega)\subset BV(\Omega)$ as in \eqref{def_bvstar}.

\begin{lemma}\label{Lem_bv_decom}
Let $\Omega$ be an open subset of $\mathbb{R}^{d}$. Then $\eta\in
BV^{*}(\Omega)$ if and only if $\eta=\eta^{(1)}+\eta^{(2)}$ where
$\eta^{(1)}\in BV(\Omega)$ with $\eta^{(1)}|_{N_{\delta '}}=0$ for some $\delta '>0$ (recall that $N_{\delta '}$ is defined in equation~\eqref{equ_Nd}), and $\eta^{(2)}\in
W^{1,2}(\Omega)$.
\end{lemma}
{\bf Proof:} We only prove the
sufficiency. Suppose $\eta\in BV^*(\Omega)$, then there exists a $\delta>0$ such that $\eta|_{N_{\delta}}\in W^{1,2}$. Let $\delta '=\delta/2$, then $\overline{N_{\delta'}}\subset N_{\delta}$. Hence there exists a Bump function $B(x)\in C^{\infty}_{c}(\Omega)$ such that $B|_{\overline{N_{\delta '}}} \equiv 1$ and $B|_{N^c_{\delta}} \equiv 0$ (here $N^c_{\delta}$ denotes the complementary of $N_{\delta}$).
We define
$$
\eta^{(2)}:=\eta\cdot B,~~\eta^{(1)}:=\eta-\eta^{(2)}.
$$
It is clear the $\eta^{(2)}\in W^{1,2}(\Omega)$ and $\eta^{(1)}\in BV(\Omega)$ with $\eta^{(1)}|_{N_{\delta '}}=0.$ \hfill $\Box$

For convenience, we denote the sequence of the Birkhoff average of $\mathcal{L}_{\overline{f}}$ on a $L^1$ function $\eta$ by $\{\eta_{n}\}_{n=1}^{\infty}$, i.e.,
$$
\eta_{n}:=\frac{1}{n}\sum_{i=0}^{n-1}\mathcal{L}_{\overline{f}}^{i}\eta.
$$
\begin{lemma}\label{lem_sobolev}
Suppose $f:X\rightarrow X$ is an invertible PWI and $\Omega\supset X$ an
open ball. Then for any $\eta\in W^{1,2}(\Omega),$ there exists a
subsequence $\{n_k\}_{k=1}^{\infty}$ such that $\eta_{n_k}\rightarrow\overline{\eta}\in
W^{1,2}(\Omega)$ in $L^{1}$ as $k\to\infty$.
\end{lemma}
{\bf Proof:} Since $\overline{f}:\Omega\to\Omega$ is a PWI, it
follows that
$\overline{f}_{i}:=\overline{f}|_{\operatorname{int}\omega_{i}}=A_{i}x+c_{i},$
where $A_{i}$ is an orthogonal matrix and $c_i$ is a translation vector. Hence, by using the orthogonality and the definition of $||\cdot||_{W^{1,2}},$
we have:
\[
||\eta_{2}\circ\overline{f}^{-i}||_{W^{1,2}}=||\eta_{2}||_{W^{1,2}}.
\]
Moreover, since $W^{1,2}(\Omega)$ is a Hilbert space, the lemma can be shown directly by using the Banach-Saks Theorem \cite{Jost05} (see Appendix~\ref{app_sobo}). \hfill $\Box$

\begin{lemma}[Coordinate Transformation]\cite{Driver03}\label{lem_divergence}
Let $\psi:W\to\Omega$ be a $C^{2}-$ diffeomorphism where $W$ and
$\Omega$ are open subsets of $\mathbb{R}^{d}.$ Given
$\overrightarrow{\phi}\in C^{1}(\Omega,\mathbb{R}^{d}),$ let
$$\overrightarrow{\phi}^{\psi}(y):=D\psi^{-1}(y)\overrightarrow{\phi}(\psi(y))$$
then
\begin{equation}\label{equ_divergence}
\operatorname{div}(|\operatorname{det}D\psi|\overrightarrow{\phi}^{\psi})=(\operatorname{div}\overrightarrow{\phi})\circ
\psi\cdot|\operatorname{det}D\psi|.
\end{equation}
\end{lemma}

\begin{lemma}\label{lem_multi_var}
Let $f:X\rightarrow X$ be an invertible
piecewise isometry with a topological partition
$\mathcal{P}=\{\omega_0,\omega_1,\cdots,\omega_{r-1}\}$. Suppose
that $\eta\in BV(\Omega)$ with $\eta|_{N_{\delta}}=0$ for some $\delta>0$, then
$\operatorname{var}(\mathcal{L}_{\overline{f}}\eta)\leq\operatorname{var}(\eta)<\infty.$
\end{lemma}
{\bf Proof:} Let $\omega_r:=\Omega\backslash X$ and
$\overline{f}_i:=\overline{f}|_{\omega_i}$ for each $\omega_{i},
i=0,\cdots,r$. For any $\eta\in BV(\Omega)$, we have $\mathcal{L}_{\overline{f}}\eta\in BV(\Omega)$. Recall that
\[
\operatorname{var}(\mathcal{L}_{\overline{f}}\eta)=\sup\left\{\int_{\Omega}(\eta\circ
\overline{f}^{-1}\cdot\operatorname{div}\overrightarrow{\phi})dm~:~\overrightarrow{\phi}\in
C_{c}^{1}(\Omega,\mathbb{R}^d),||\overrightarrow{\phi}||_{\infty}\leq
1\right\}.
\]
Hence, given any $\epsilon>0$, there exists $\overrightarrow{\phi}\in
C_{c}^{1}(\Omega,\mathbb{R}^d)$ with $||\overrightarrow{\phi}||_{\infty}\leq 1$
such that
\begin{eqnarray}\label{equation directly}
\operatorname{var}(\mathcal{L}_{\overline{f}}\eta)-\epsilon\leq\int_{\Omega}(\eta\circ
\overline{f}^{-1}\cdot \operatorname{div}\overrightarrow{\phi})dm
\nonumber
=\sum_{i=0}^{r}\int_{\overline{f}_{i}(\operatorname{int}\omega_{i})}(\eta\circ\overline{f}_{i}^{-1}\cdot\operatorname{div}\overrightarrow{\phi})dm.
\end{eqnarray}
Using coordinate transformation $(x=\overline{f}_{i}^{-1}(y))$ on each $\omega_{i}$, we have
\[
\int_{\overline{f}_{i}(\operatorname{int}\omega_{i})}\eta(\overline{f}_i^{-1}(y))\cdot\operatorname{div}\overrightarrow{\phi}(y)dm(y)
=\int_{\operatorname{int}\omega_{i}}\eta(x)\cdot\operatorname{div}(\overrightarrow{\phi})(\overline{f}_i(x))dm(\overline{f}_i(x))
\]
\[
=\int_{\operatorname{int}\omega_i}\eta\cdot(\operatorname{div}\overrightarrow{\phi})\circ
\overline{f}_i
dm=\int_{\operatorname{int}\omega_{i}}\eta\cdot\operatorname{div}\overrightarrow{\phi}^{\overline{f}_{i}}dm,
\]
where
$\overrightarrow{\phi}^{\overline{f}_{i}}:=(D\overline{f}_{i})^{-1}\cdot\overrightarrow{\phi}\circ
\overline{f}_{i}$. The second equality is due to
$|\operatorname{det}D\overline{f}_{i}(x)|\equiv 1$ while the third is due to equation~\eqref{equ_divergence} in Lemma~\ref{lem_divergence}. Hence,
$$
\operatorname{var}(\mathcal{L}_{\overline{f}}\eta)-\epsilon\leq\sum_{i=1}^{r}\int_{\operatorname{int}\omega_i}\eta(x)\cdot\operatorname{div} \overrightarrow{\phi}^{\overline{f}_{i}}dm.
$$

Moreover, each $\overline{f}_i$ can be written in the form of $\overline{f}_i(x)=A_i\cdot x+c_i$ where the orthogonal matrix $A_i$ preserves the Euclidean metric. Then
$$
\sup_{x\in\operatorname{int}\omega_{i}}|\overrightarrow{\phi}^{\overline{f}_{i}}(x)|=\sup_{x\in\operatorname{int}\omega_{i}}|A_i^{-1}\overrightarrow{\phi}(\overline{f}_{i}(x))|=
\sup_{x\in\operatorname{int}\omega_{i}}|\overrightarrow{\phi}(\overline{f}_{i}(x))|
=\sup_{y\in\overline{f}_{i}(\operatorname{int}\omega_{i})}|\overrightarrow{\phi}(y)|\leq
1.
$$

Since $\eta|_{N_{\delta}}=0$ for some $\delta>0$, for each $\omega_i$, we let $\omega_i^{(\delta)}:=N^c_{\delta}\cap\omega_i$ which is a compact subset of $\omega_i$. Moreover, $\omega_{i}^{(\delta)}\subset\omega_{i}^{(\delta/2)}\subset\operatorname{int}\omega_{i}$. Therefore, there exists a bump function $B_{i}(x)\in
C_{c}^{\infty}(\Omega)$ such that $0\leq
B_{i}(x)\leq1$, $B_{i}(x)=1$ on
$\omega_{i}^{(\delta)}$ and $B_{i}(x)=0$ on
$(\operatorname{int}\omega^{(\delta/2)})^{c}$. We extend the function $\overrightarrow{\phi}^{\overline{f}_i}$ to be zero outside $\omega_i$ and define
\begin{equation*}
   \overrightarrow{\psi}(x):=\sum_{i=0}^r B_{i}(x)\cdot\overrightarrow{\phi}^{\overline{f}_i}(x), ~x\in\Omega.
\end{equation*}
It is apparent to see that $\overrightarrow{\psi}\in
C_{c}^{1}(\Omega,\mathbb{R}^{d})$ and
$\sup_{x\in\Omega}|\overrightarrow{\psi}|\leq\sup_{x\in\Omega}|\overrightarrow{\phi}^{\overline{f}_i}(x)|\leq 1$ for each $i$. Moreover, since $\overrightarrow{\psi}|_{\omega_{\delta}}=\overrightarrow{\phi}^{\overline{f}_i}|_{\omega_{\delta}}$ and $\eta|_{N_{\delta}}=0$, we have
$$
\operatorname{var}(\mathcal{L}_{\overline{f}}\eta)-\epsilon\leq
\sum_{i=1}^{r}\int_{\operatorname{int}\omega_i}\eta(x)\cdot\operatorname{div}
\overrightarrow{\phi}^{\overline{f}_{i}}dm=\int_{\Omega}\eta\cdot\operatorname{div}\overrightarrow{\psi}dm\leq\operatorname{var}(\eta).
$$
Since $\epsilon$ is arbitrary, it follows that
$\operatorname{var}(\mathcal{L}_{\overline{f}}\eta)\leq\operatorname{var}(\eta)<\infty$. \hfill $\Box$

{\bf \ \ Proof of Theorem \ref{theo_2}:} Given any $\eta\in
BV^*(\Omega)$, from Lemma~\ref{Lem_bv_decom}, we can write
$\eta=\eta^{(1)}+\eta^{(2)}$ where $\eta^{(1)}|_{N_{\delta}}=0$ for
some $\delta>0$ and $\eta^{(2)}\in W^{1,2}$. Moreover, by
Lemma~\ref{lem_sobolev}, there exists a subsequence
$\{\eta^{(2)}_{n_k}\}$ of $\{\eta^{(2)}_{n}\}$ which converges to a
function $\overline{\eta}^{(2)}\in W^{1,2}$ in $L^1(m)$.

Consider the function $\eta^{(1)}$. Given any $n\geq 1$, by analogous arguments used in Lemma~\ref{lem_multi_var}, we have
$\operatorname{var}(\mathcal{L}_{\overline{f}}^{n}\eta^{(1)})\leq\operatorname{var}(\eta^{(1)})$.
It follows that $||\mathcal{L}_{\overline{f}}^n\eta^{(1)}||_{BV}\leq||\eta^{(1)}||_{BV}$. Hence,
$||\eta^{(1)}_{n}||_{BV}\leq||\eta^{(1)}||_{BV}<\infty$.
By Helly's Theorem~\cite{Evans92}, there exists a subsequence, for convenience say,
$\{\eta^{(1)}_{n_{k}}\}$ converging to a function
$\overline{\eta}^{(1)}\in BV$ in $L^{1}(m)$.

By the triangle inequality, for each $\overline{\eta}^{(i)}$, $i=1,2$,
\begin{equation}\label{equ_triangle}
||\mathcal{L}_{\overline{f}}\overline{\eta}^{(i)}-\overline{\eta}^{(i)}||_{1}\leq
||\mathcal{L}_{\overline{f}}\overline{\eta}^{(i)}-\mathcal{L}_{\overline{f}}\eta^{(i)}_{n_{k}}||_{1}+||\mathcal{L}_{\overline{f}}\eta^{(i)}_{n_{k}}-\eta^{(i)}_{n_{k}}||_{1}
+||\eta^{(i)}_{n_{k}}-\overline{\eta}^{(i)}||_{1}.
\end{equation}
It is clear that the first and third term tend to be $0$ as $k\to\infty$.
Moreover, we note that
\[
||\mathcal{L}_{\overline{f}}\eta^{(i)}_{n_{k}}-\eta^{(i)}_{n_{k}}||_{1}=\left|\left|\frac{1}{n_{k}}\sum_{j=1}^{n_{k}}\eta^{(i)}\circ\overline{f}^{-j}-\frac{1}{n_{k}}\sum_{j=0}^{n_{k}-1}\eta^{(i)}\circ
\overline{f}^{-j}\right|\right|_{1} =\frac{1}{n_{k}}||\eta^{(i)}\circ
\overline{f}^{-n_{k}}-\eta^{(i)}||_{1}\leq\frac{2}{n_{k}}||\eta^{(i)}||_{1},
\]
hence the second term tends to $0$ as $k\rightarrow \infty$. Therefore,
$\mathcal{L}_{\overline{f}}\overline{\eta}^{(i)}=\overline{\eta}^{(i)}$,
implying that $\overline{\eta}:=\overline{\eta}^{(1)}+\overline{\eta}^{(2)}\in
BV(\Omega)$ in an invariant density of $\overline{f}$ from Lemma~\ref{lemma_fixed}. Moreover, since $\eta|_{X}\geq 0$ and $||\overline{\eta}|_X||_1=||\eta|_X||_1>0$, by normalization $d\mu=\overline{\eta}|_{X}dm\in \mathcal{M}_{IB}(f).$ \hfill$\Box$

\begin{remark}
For a given invertible PWI, if the set $N_{\delta}$ satisfies
$f^{-1}(N_{\delta})\subset N_{\delta}$, then the above accumulation
point $\overline{\eta}$ can show to be in $BV^*(\Omega)$ for the
following reason. Since the sequence $\{\eta^{(1)}_{n_{k}}\}$
converges to $\overline{\eta}^{(1)}$ in $L^{1}$, then there will
exist a subsequence which pointwise converges to
$\overline{\eta}^{(1)}. $ Moreover, since
$\eta^{(1)}|_{N_{\delta}}=0$ and
$\overline{f}^{-1}(N_{\delta})\subset N_{\delta}$, it follows that
$\eta^{(1)}_n|_{N_{\delta}}=0$ for any $n\geq 0$. Therefore,
$\overline{\eta}^{(1)}|_{N_{\delta}}=0.$ By Lemma
\ref{Lem_bv_decom}, this implies
$\overline{\eta}=\overline{\eta}^{(1)}+\overline{\eta}^{(2)}\in
BV^*(\Omega)$.
\end{remark}

\subsection{Proof of Proposition~\ref{prop_density}}
\begin{lemma}\label{lem_Condition_expectation}
Let $(X, \mathfrak{B},m)$ be a probability space and $f:X\rightarrow
X$ be an invertible PAP, then $d\mu:=\varphi dm\in
\mathcal{M}_{I}(f)$ if and only if
$\varphi=\mathbb{E}(\varphi|\mathcal{I})$, where $\mathbb{E}(\varphi|\mathcal{I})$ is the conditional expectation on the $\sigma-$field $\mathcal{I}:=\{B\in\mathfrak{B}|f^{-1}(B)=B~\operatorname{mod}~m\}$.
\end{lemma}

{\bf Proof:} Suppose $d\mu=\varphi dm\in\mathcal{M}_{I}(f)$, then by Lemma~\ref{lemma_fixed} we have $\mathcal{L}_{f}\varphi=\varphi$, which implies $\varphi\circ f=\varphi$. Combining with the Birkhoff Ergodic Theorem, we obtain
\begin{equation*}
\varphi=\lim_{n\rightarrow\infty}\frac{1}{n}\sum_{i=0}^{n-1}\varphi\circ
f^{i}=\mathbb{E}(\varphi|\mathcal{I}).
\end{equation*}

For the converse, given any fixed
$L\in\mathbb{N},$ we define $\varphi_{L}(x):=\min\{\varphi(x),L\}$.
By the Birkhoff Ergodic Theorem, for $m-a.e.~x\in X$,
\begin{eqnarray*}
\mathcal{L}_{f}(\mathbb{E}(\varphi_{L}|\mathcal{I}))(x)=\mathbb{E}(\varphi_{L}|\mathcal{I})\circ
f^{-1}(x)=\lim_{n\rightarrow\infty}\frac{1}{n}\sum_{i=-1}^{n-2}\varphi_{L}\circ
f^{i}(x)\\
=\lim_{n\rightarrow\infty}\left[\frac{1}{n}(\sum_{i=0}^{n-1}\varphi_{L}\circ
f^{i})(x)+\frac{1}{n}\left(\varphi_{L}\circ
f^{-1}(x)-\varphi_{L}\circ f^{n-1}(x)\right)\right].
\end{eqnarray*}
Since
$$
\lim_{n\rightarrow\infty}\frac{1}{n}\left(\varphi_{L}\circ
f^{-1}(x)-\varphi_{L}\circ f^{n-1}(x)
\right)\leq\lim_{n\to\infty}\frac{2L}{n}=0,
$$
we have $\mathcal{L}_{f}\mathbb{E}(\varphi_{L}|\mathcal{I})=\mathbb{E}(\varphi_{L}|\mathcal{I}).$
Then by the Monotone convergence theorem,
\begin{equation*}
\mathcal{L}_{f}\mathbb{E}(\varphi|\mathcal{I})=\mathbb{E}(\varphi|\mathcal{I}).
\end{equation*}
Therefore, $\varphi=\mathbb{E}(\varphi|\mathcal{I})$ implies $d\mu:=\varphi dm\in\mathcal{M}_I(f)$ from Lemma~\ref{lemma_fixed}.
\hfill $\Box$

We say a $\sigma-$field $\mathcal{I}$ is {\em finitely generated}, if there exists a partition
$\mathcal{A}:=\{A_{i}\}_{i=0}^{r-1}\subseteq\mathcal{I}$,
and for each $B\in\mathcal{I},$ there exist finitely many
$A_{i_{1}},\cdots,A_{i_{l}}\in\mathcal{A}$ such that
$
B=\bigcup_{k=1}^{l}A_{i_{k}}\operatorname{mod}m.
$

{\bf Proof of Proposition \ref{prop_density}:} We first show statement $(i)$. Since there are only finitely many ergodic measures
$\{\nu_{i}\}_{i=0}^{r-1}\in \mathcal{M}_{I}(f)$ for any topologically transitive IET $f$ \cite{Keane75}, the $\sigma-$field $\mathcal{I}$ is finitely generated by a partition, say $\mathcal{A}:=\{A_{0},\cdots,A_{r-1}\}$ with $0<m(A_i)<1$ and $\nu_i=m|_{A_i}\in\mathcal{M}_I(f)$ is ergodic for each $i$. For any $d\mu=\varphi dm\in\mathcal{M}_{I}$, by Lemma~\ref{lem_Condition_expectation} and \cite{Boy97},
\begin{equation}\label{equ_density}
\varphi=\mathbb{E}(\varphi|\mathcal{I})=\sum_{i=0}^{r-1}\frac{1}{m(A_{i})}\int_{A_{i}}\varphi
dm\cdot \chi_{A_{i}},~~\forall~ \varphi\in L^{1}(m).
\end{equation}
Therefore, $\varphi$ is a simple function.

Now we show statement $(ii)$. Since each $m|_{A_i}\in\mathcal{M}_{I}(f)$ and $m(A_{i})>0$, it is clear that for $m-a.e.~x\in A_i$, the orbit $O_f(x)\subset A_i$. Moreover, since $f$ is minimal, each $A_{i}$ is dense in $X$. Hence $\operatorname{int}A_{i}=\emptyset.$ Based on
\eqref{equ_density}, if $\mu:=\varphi dm\in \mathcal{M}_{I}(f)$ and
$\mu\neq m$ then there exists $i\neq j$ such that
$\int_{A_{i}}\varphi
dm/m(A_{i})\neq\int_{A_{j}}\varphi dm/m(A_{j}).$ Therefore,
$\varphi$ is discontinuous everywhere. Moreover, for any open set $U\subset [0,1)$,  we have $m(A_i\cap U)>0$. Hence $\mu(U)>0$. This follows that $\operatorname{supp}(\mu)=[0,1)$. \hfill$\Box$

\subsection{Proofs of Lemma~\ref{lemma_piece_invar} and Proposition~\ref{pro_nonin}}\label{sec_pro_noninver}

{\bf Proof of Lemma~\ref{lemma_piece_invar}:} We first show almost closedness of $X^+$. Since $f_{i}$ is Lipschitz continuous, it
can be continuously extended from $\operatorname{int}\omega_{i}$
onto $\overline{\operatorname{int}\omega_{i}}$. If we denote its
continuous extension by
$\widehat{f}_{i}:\overline{\operatorname{int}\omega_{i}}\to
\widehat{f}_{i}(\overline{\operatorname{int}\omega_{i}})$, then each $\widehat{f}_{i}$ is also Lipschitz continuous. Together with the non-singularity of $f$, we know \footnote{By $A\subset B$ mod $m$ we mean that for $m-a.e.$ $x\in A$, we have $x\in B$.}
\begin{equation*}
\begin{array}{lll}
\overline{X^{+}}:=\operatorname{closure}\left(\bigcap_{j=0}^{\infty}f^{j}(\bigcup_{i=0}^{r-1}\omega_{i})\right)
&\subseteq&\bigcap_{j=0}^{\infty}\operatorname{closure}\left(\bigcup_{i=0}^{r-1}\widehat{f}_{i}^{j}(\operatorname{int}\omega_{i})\right)~\mbox{mod}~m\\
& \subseteq&\bigcap_{j=0}^{\infty}\bigcup_{i=0}^{r-1}\widehat{f}_{i}^{j}(\overline{\operatorname{int}\omega_{i}})~\mbox{mod}~m \\
&= & \bigcap_{j=0}^{\infty}\bigcup_{i=0}^{r-1}f_i(\omega_i) \mbox{~mod~}m\\
&=&X^+ \mbox{~mod~} m.
\end{array}
\end{equation*}
This implies $m(X^{+})=m(\overline{X^{+}}).$

Denote $\omega^{+}_{i}:=\omega_{i}\cap X^{+}$, then $X^{+}=\cup_{i=0}^{r-1} \omega^{+}_i$. Consequently,
$$
\sum_{i=0}^{r-1}m(f(\omega^{+}_{i}))=\sum_{i=0}^{r-1}m(\omega^{+}_{i})=m(\bigcup_{i=0}
^{r-1}\omega^{+}_{i})=m(X^+)=m(f(X^{+}))=m(\cup_{i=0} ^{r-1}
f(\omega^{+}_{i})).
$$
This follows that $m(f(\omega_{i}^{+})\cap
f(\omega_{j}^{+}))=0,$ implying $f|_{X^+}$ is $m-a.e.$
invertible. \hfill$\Box$

{\bf Proof of Proposition~\ref{pro_nonin}:} We note that $f^+$ is piecewise Lipschitz continuous, then for any Borel subset $A\subset \overline{X^{+}}$ with $m(A)=0$, we have $m(f^{+}(A))=0$. Moreover, since $f|_{\omega_i}$ is bi-Lipschitz, it follows $m((f^{+})^{-1}(A))=0$. Hence $f^+$ is non-singular. The $m-a.e.$ invertibility of $f^+$ is directly from Lemma~\ref{lem_almostclosedness}.

When $m(\overline{X^{+}})>0$, $f|_{\overline{X^{+}}}$ can be
viewed as a first return map of $f$ on $\overline{X^{+}}$. Therefore,
$$\left\{\mu(\cdot):=\nu(\cdot\cap\overline{X^+}),\forall\nu\in\mathcal{M}_{I}(f^+)\right\}\subseteq\mathcal{M}_{I}(f).$$
Moreover, for any $\mu\in\mathcal{M}_{I}(f),$ since
$\mu(X)=1$ and $f^{-1}\circ f(X)=X$, hence $\mu(f(X))=1$. This
implies that $\mu(X^{+})=1$ and completes the proof of statement
(i).

For statement (ii), we prove by contradiction. Suppose that there
exists $\mu\in \mathcal{M}_{I}(f)$, since $m(\overline{X^{+}})=0$, it follows that $\mu(\overline{X^{+}})=0$. This is a
contradiction with $\mu(X^{+})=1$. \hfill $\Box$

\section*{Acknowledgments}
The authors are grateful to Prof. Peter Ashwin, Prof. Henk Bruin,
Dr. Mark Holland and Dr. Corinna Ulcigrai as well as an anonymous referee for their very helpful comments. The authors also thank the symposium on Ergodic Theorem and Dynamical Systems organized in University of Warwick during 2010-2011.
\appendix
\section{Appendix}

\subsection{Bounded variation} \label{app_BV}

We introduce the usual notion of bounded variation in one dimension
followed by definitions of multidimensional bounded variation.

Let $\eta\in L^1(\mathbb{R})$ and $[a,b]$ be an interval outside of which $\eta(x)=0$. The {\em total variation} of $\eta$ is defined to be
\[
V(\eta):=\sup\sum_{i=0}^{i-1}|\eta(x_{i+1})-\eta(x_{i})|,
\]
where ``$\sup$'' is taken over all possible finite partitions of the interval $[a,b]$ by points $x_{0}=a<x_{1}<\cdots<x_{r}=b$. The {\em essential total variation} of $\eta$ is defined as
$
\bar{V}(\eta):=\inf_u V(\eta+u),
$
where the ``$\inf$'' is taken over all functions $u$ that equal zero
almost everywhere on $[a,b]$ \cite{Volper85}.

Next, we proceed to describe a definition of multidimensional
bounded variation. Let $\Omega$ be an open subset of $\mathbb{R}^d$
and $\eta\in L^1(\Omega)$ be a function with compact support. We regard
$\eta(x)$ as a function of the variable $x_{i}$ for the other
variables fixed and denote by $\bar{V}_{i}(x_{i}')$ the essential
total variation of the function $\eta$ with respect to $x_{i}$ for a
fixed point
$x_{i}'=\{x_{1},x_{2},\cdots,x_{i-1},x_{i+1},\cdots,x_{d}\}.$
\begin{definition}\cite{Volper85}\label{def_induc}
Suppose $\Omega$ is an open subset of $\mathbb{R}^d$. A function $\eta\in L^1(\Omega)$ with compact support is said to be of {\em bounded variation} if the integrals
$$
\int \bar{V}_i(x_i')dx_i'<\infty, ~~~~(i=1,2,\cdots,n).
$$
\end{definition}

We give a further definition of bounded variation below.
\begin{definition}\cite[p158]{Volper85}\label{def_functional} Suppose $\Omega$ is an open subset of $\mathbb{R}^{d}$ and a function $\eta\in L^{1}(\Omega),$ then $\eta$ is said to be of {\em bounded variation} if there exists a constant $K$ such that
\begin{equation}\label{equ_functional}
  \left |\int_{\Omega}\frac{\partial\phi}{\partial x_{i}}\eta dx\right|\leq K\sup_{x\in
\Omega}|\phi(x)|,~~~(i=1,2,\cdots,d)
\end{equation}
for all $\phi\in C^{1}_{c}(\Omega,\mathbb{R})$.
\end{definition}

Definition~\ref{def_induc} and Definition~\ref{def_functional} coincide for functions $\eta\in L^1(\Omega)$ with compact support \cite{Volper85}. Moreover, we show the equivalence between Definition~\ref{def_functional} and Definition~\ref{def_bv} in Section~\ref{sec_bv} via the following proposition.

\begin{proposition}\label{prop_simple equivent}
Suppose $\Omega$ is an open subset of $\mathbb{R}^n$ and $\eta\in L^{1}(\Omega)$. Then
$\operatorname{var}(\eta)<\infty$ (recall that
$\operatorname{var}(\cdot)$ is defined in Definition~\ref{def_bv}) if
and only if the inequality~\eqref{equ_functional} holds for all $\phi\in
C^{1}_{c}(\Omega,\mathbb{R})$.
\end{proposition}
{\bf Proof:}~~``$\Rightarrow$''. Suppose that
$\operatorname{var}(\eta)<\infty$, then by Definition~\ref{def_bv},
\begin{equation}\label{equ_bounded}
 \left|\int_{\Omega}\eta(x)\operatorname{div}\overrightarrow{\phi}(x)dx\right|\leq
\operatorname{var}(\eta)||\overrightarrow{\phi}||_{\infty},~~\forall\overrightarrow{\phi}\in
C_{c}^{1}(\Omega,\mathbb{R}^{d}),
\end{equation}
For any $\phi\in C_{c}^{1}(\Omega,\mathbb{R})$, let $\overrightarrow{\phi_{i}}=(\overbrace{0,\cdots}^{i-1},\phi,\cdots,0)$, then the inequality~\eqref{equ_bounded} holds for $\overrightarrow{\phi}_i$. This implies the inequality~\eqref{equ_functional} holds for all $C_c^1(\Omega,\mathbb{R}^d)$.

``$\Leftarrow$''. By the inequality~\eqref{equ_functional} and Definition~\ref{def_bv}, it is clear that
$\operatorname{var}(\eta)\leq d\cdot K<\infty$ where $d$ is the dimension
constant. \hfill $\Box$

\subsection{Sobolev space $W^{1,2}$}\label{app_sobo}
\begin{definition}\cite{Jost05}
Suppose $\Omega$ is an open set in $\mathbb{R}^d$. The Sobolev space $W^{1,2}(\Omega)$ is defined to be the set of all functions $u\in  L^2(\Omega)$ such that for every multi-index $\alpha=(\alpha_1,\cdots,\alpha_d)$ with $\alpha_i\geq 0$ and $|\alpha|:=\sum_{j=1}^{d}\alpha_{j}\leq 1$, the weak partial derivative $D_{\alpha}u$ belongs to $L^2(\Omega)$, i.e.,
$$
W^{1,2}(\Omega) = \{ u \in L^2(\Omega) : D_{\alpha}u \in L^2(\Omega), \forall |\alpha| \leq 1 \}.
$$
We define a norm on $W^{1,2}(\Omega)$ by
$$
||u||_{W^{1,2}}:=(\sum_{|\alpha|\leq
1}\int_{\Omega}|D_{\alpha}u|^{2}dm)^{1/2}.
$$
\end{definition}

The Sobolev space $W^{1,2}(\Omega)$ with norm $||\cdot||_{W^{1,2}}$ is a Hilbert subspace of $BV(\Omega)$. We state the following Banach-Saks theorem which is applied to the sobolev space $W^{1,2}$ in Lemma~\ref{lem_sobolev}.

{\bf Banach-Saks Theorem}~\cite{Jost05}
 Let $\{x_n\}_{n\mathbb{N}}$ be a sequence from a Hilbert space with $||x_n||\leq K$ (independent of $n$), then there exists a subsequence $\{x_{n_j}\}_{j\in\mathbb{N}}$ and an $x\in H$ such that
$$
\frac{1}{k}\sum_{j=1}^{k}x_{n_j}\rightarrow x,~~\mbox{as~}k\rightarrow \infty
$$
in norm.

\subsection{Piecewise rotations}\label{app_piecerotation}
\begin{definition}\cite{Goetz01}
Let $X$ be a compact subset of $\mathbb{C}$. A map $T: X\rightarrow X$ is called a {\em piecewise
rotation} with a partition
$\mathcal{P}:=\{\omega_{0},\cdots,\omega_{r-1}\}$ if
\[
T|_{\omega_j}x=\rho_{j}x+z_{j},~~~x\in\omega_{j}
\]
for some complex numbers $z_{j}$ and $\rho_{j}$ such that
$|\rho_{j}|=1$ for all $j=0,1,\cdots,r-1$. The atoms are assumed to
be mutually disjoint convex polygons.
\end{definition}

It is clear that piecewise rotations are PWIs in
$\mathbb{R}^{2}$ with a topological partition $\mathcal{P}$ and are homeomorphisms when restricting on each atom.

\subsection{Interval translation maps and interval exchange transformations}\label{app_IET}
\begin{definition}
Let $I=[0,1)$ be an interval and $0=\beta_{0}<\beta_{1}<\cdots<\beta_{r}=1$ be a finite partition of $I$. An interval map $T:I\rightarrow I$ is said to be an {\em interval translation map} \cite{Bruin03} if
\[
T(x)=x+\gamma_{i},~~\beta_{i-1}\leq x<\beta_i,
\]
where each $\gamma_{i}$ is a fixed real number.
Particularly, if $T$ maps $I$ onto itself then $T$ is called an {\em interval exchange
transformation.}
\end{definition}

\end{document}